\documentclass[12pt]{article}

\newtheorem{thm}{Theorem}
\newtheorem{lem}[thm]{Lemma}

\def\b/{$B_1^n$}
\def\bb/{$B_2^n$}
\def\elip/{ellipsoid}

\def\C{{\bf C}}
\def\BM/{Brunn-Minkowski}
\def\bod/{convex body}
\def\bods/{convex bodies}
\def\d/{Dvoretzky}
\def\ep/{\varepsilon}
\def\iq/{inequality}
\def\iqs/{inequalities}
\def\iso/{isoperimetric}

\def\P/{probability}
\def\sph/{$S^{n-1}$}

\def\startproof{\bigskip \noindent {\em Proof }}
\def\endproof{\hfill $\vcenter{\hrule height .3mm
		\hbox {\vrule width .3mm height 2mm \kern 2mm
			\vrule width .3mm} \hrule height .3mm}$ \bigskip}

\title{Rational approximations to the zeta function}
\author{Keith Ball}

\begin{document}

	\maketitle
	
	\begin{abstract}
		This article describes 
		a sequence of rational functions which
		converges locally uniformly to $\zeta$. The numerators (and denominators)
		of these rational functions can be expressed as characteristic polynomials of matrices that are on the face of it very simple.
		As a consequence, the Riemann hypothesis can be restated as what looks like a rather conventional spectral problem but which is related to the one found by Connes and by Berry and Keating. However the point here is that the rational approximations look to be susceptible of quantitative estimation.
	\end{abstract}
	
	\section*{Introduction}

	This article describes 
	a sequence of rational functions that converges locally uniformly to $\zeta$ at least to the right of the line 	$\{s: \Re s =0 \}$. The sequence begins
	\[ \frac{1}{(s-1)}, \; \; \;  \frac{s+1}{2(s-1)}, \; \; \; \frac{4 s^2 +11 s+9}{6 (s+3)(s-1)}, \; \; \;
	\frac{(s+2)(3 s^2+10 s+11)}{4 (s^2+6 s+11)(s-1)} ,\]
	\begin{equation}
	\label{funs} 
	\frac{(s+2)(72 s^3+490 s^2+1193 s+1125)}{30 (3 s^3+29 s^2+106 s+150)(s-1)}, \ldots
	\end{equation}
	
	We define the sequence as follows. For each integer $m \geq 0$ we define 
	\[ p_m(t)= \left(1-t \right) \left(1-\frac{t}{2} \right) \ldots \left( 1-\frac{t}{m} \right) \]
	and the coefficients $(a_{m,j})$ by
	\[ p_m(t)=\sum_0^m (-1)^j a_{m,j} t^j. \]
	These coefficients are rescaled Stirling numbers of the first kind but it is more convenient for us to use the different indexation and scaling here.
	We then set
	\[ F_m(s)=\sum_0^m \frac{a_{m,j} B_j}{s+j-1} \]
	where the $B_j$   are the usual Bernoulli numbers given by the generating function
	\[ \frac{z}{e^z-1} = \sum_0^\infty \frac{B_j}{j!} z^j \]
	and
	\[ G_m(s)=\sum_{j=0}^m (-1)^j \frac{a_{m,j}}{s+j-1}. \]
	The rational functions in question are the ratios
	\[ \frac{F_m(s)}{(s-1)G_m(s)}. \]
	For example
	\[ F_3(s)= \frac{1}{s-1} - \frac{11}{12 s} + \frac{1}{6 (s + 1)} = 
	\frac{3s^2+10s+11}{12 (s-1)s(s+1)} \]
	and 
		\[ G_3(s)= \frac{1}{s-1} - \frac{11}{6 s} + \frac{1}{s + 1}-\frac{1}{6(s+2)} = 
		\frac{s^2+6s+11}{3 (s-1)s(s+1)(s+2)}. \]
		
	The main theorem of the article is as follows.
	\begin{thm}
		\label{main}
	For each $m$ set
	 \[ F_m(s)=\sum_0^m \frac{a_{m,j} B_j}{s+j-1} \]
			where the $B_j$   are the usual Bernoulli numbers 
			and
			\[ G_m(s)=\sum_{j=0}^m (-1)^j \frac{a_{m,j}}{s+j-1} \]
	where
	\[ \left(1-t \right) \left(1-\frac{t}{2} \right) \ldots \left( 1-\frac{t}{m} \right) 
	=\sum_0^m (-1)^j a_{m,j} t^j. \]
	Then
	\[ \frac{F_m(s)}{(s-1)G_m(s)} \rightarrow \zeta(s) \]
	locally uniformly on $\{s: \Re s >0 \}$ with the obvious convention at $s=1$.
		\end{thm}

	It is immediate that the $m^{th}$ ratio interpolates $\zeta$ at the points $0,-1,-2,\ldots,1-m$ and has a simple pole with residue 1 at $s=1$. We shall show that the numerators (and denominators)
	of these rational functions can be expressed as characteristic polynomials of matrices that are on the face of it very simple.
	One way to state this is that the numerator of the $m^{th}$ function is the determinant of
	\[ \left( \begin{array}{cccccc} 1 & 0 & 0 & 0 & \ldots & 0 \\
	\frac{1}{2} & 1 & 0 & 0 & \ldots & 0 \\
	\frac{1}{3} & \frac{1}{2} & 1 & 0 & \ldots & 0 \\
	\frac{1}{4} & \frac{1}{3} & \frac{1}{2} & 1 &  & 0 \\
	\vdots & & & \ddots & \ddots & \vdots \\
	\frac{1}{m+1} & \frac{1}{m} & \frac{1}{m-1} & \ldots & \frac{1}{2} & 1 \end{array} \right) +
	(1-s) \left( \begin{array}{cccccc} 0 & 1 & \frac{1}{2} & \frac{1}{3} & \ldots & \frac{1}{m} \\
	0 & 0 & \frac{1}{2} & \frac{1}{3} & \ldots & \frac{1}{m}\\
	0 & 0 &  0 & \frac{1}{3} & \ldots & \frac{1}{m} \\
	0 & 0 &  0 & 0 & & \vdots \\
	\vdots & \vdots & \vdots &  & \ddots & \frac{1}{m} \\
	0 & 0 & 0 & 0 & \ldots & 0 \end{array} \right). \]
	
	The real content of the Riemann hypothesis is that $\zeta$ has no zeroes to the {\em right} of the critical line. The customary formulation is equivalent by virtue of the functional equation proved by Riemann. In order to show that a holomorphic function has no zeroes in (say) a half-plane it suffices to express the function as a locally uniform limit of holomorphic functions with no zeroes there.
	
	As a consequence, the Riemann hypothesis can be restated as what looks like a rather conventional spectral problem: to show that the spectra of certain matrices stay to the left of the critical line, at least asymptotically as $m \rightarrow \infty$. Needless to say, I spent some time thinking about this problem without success but am {\em certain} that I have not exhausted the possible lines of attack. Even if the simplicity of these matrices is indeed an illusion, as one would expect, there are concrete reasons to think that these rational approximations to zeta might be useful for estimating the size of zeta, in the sense of the Lindel\"{o}f hypothesis: see for example the book of Patterson \cite{Pat}.
	
	Poly\'{a} suggested that the Riemann hypothesis should be proved by expressing the zeroes of zeta (rotated onto the real line) as eigenvalues of a self-adjoint operator. (The statement is often credited to Hilbert: it was Terry Tao who pointed out the mistake to me.) A number of candidates for such operators have been proposed, coming mainly from quantum theory. The best known of these were found by Connes \cite{con} and by Berry and Keating \cite{BK}. There is a connection between these infinite-dimensional operators and the finite-dimensional ones described here, which will be explained briefly in Section \ref{Connes}. This became apparent to me from the very readable article of Lachaud \cite{lach}. My hope is that the finite-dimensional operators and the rational functions they correspond to, are susceptible of quantitative estimation that would not make sense for the infinite-dimensional operators.
	
	It is well known that the zeroes of $\zeta$ should be modelled by the eigenvalues of certain random matrices. This originated in the work of Dyson and Montgomery \cite{M}
	and was experimentally confirmed by remarkable calculations of Odlyzko \cite{Odl}. Katz and Sarnak extended the model to other $L$-functions \cite{Kat}. 
	In the past two decades a huge amount of work has been done on this connection in particular by Keating, Snaith and their collaborators \cite{Keat}. 
	While strictly speaking this is only indirectly related to the results in this article the random model is clearly a crucial inspiration.
	
	In order to prove the convergence of $\frac{F_m(s)}{(s-1)G_m(s)}$ to $\zeta(s)$ we shall show that 
	\[ F_m(s) \approx h_m^{1-s} \Gamma(s) \zeta(s) \]
	and
	\[ (s-1)G_m(s) \approx h_m^{1-s} \Gamma(s)  \]
	where $h_m$ is the partial sum $\sum_{j=1}^m 1/j$ of the harmonic series. 
	We can get a sense of why this is, quite easily.
	Using the following variant of Kronecker's formula 
	\begin{equation}
	\label{Bern} 
	B_j = (-1)^j \sum_{k=0}^j \frac{1}{k+1} \sum_{r=0}^k {k \choose r} (-1)^r (r+1)^j. 
	\end{equation}
	it is easy to show that for $\Re s>1$
	\begin{equation}
	\label{Fform}
	F_m(s)= \int_0^1 \left(\sum_{k=0}^m \frac{1}{k+1} \sum_{r=0}^k {k \choose r} (-1)^r p_m((r+1) x) 
	\right) x^{s-2} \, dx. 
	\end{equation}
	
	If $x$ is close
	to zero then $p_m(x)$ is approximately 
	\[ \prod_{i=1}^m e^{-x/i} = e^{-h_m x}. \]
	So the sum over $r$ in equation (\ref{Fform}) is approximately
	\[ \sum_{r=0}^k (-1)^r {k \choose r} e^{-h_m (r+1) x} = e^{-h_m x} (1-e^{-h_m x})^k. \]
	Thus for small values of $x$ the integrand in equation (\ref{Fform}) is approximately
	\[ \left(\sum_{k=0}^m \frac{1}{k+1} e^{-h_m x} (1-e^{-h_m x})^k \right) x^{s-2}. \]
	If the approximation were good for {\em all} $x$ between $0$ and $1$ then $F_m(s)$ would be close to
	\[ \int_0^1 \sum_{k=0}^m \frac{1}{k+1} e^{-h_m x} (1-e^{-h_m x})^k  x^{s-2} \, dx = 
	h_m^{1-s} \int_0^{h_m} \sum_{k=0}^m \frac{1}{k+1} e^{-y} (1-e^{-y})^k  y^{s-2} \, dy. \]
	The last integral plainly converges to 
	\[  \int_0^\infty \sum_{k=0}^\infty 
	\frac{1}{k+1} e^{-y} (1-e^{-y})^k y^{s-2} \, dy =\int_0^\infty \frac{y}{1-e^{-y}} e^{-y} y^{s-2} \, dy \] as $m \rightarrow \infty$ provided
	$\Re s>1$ and the latter is easily seen to be $\Gamma(s) \zeta(s)$. 
	
	Our first aim will be to show that indeed
	\[ h_m^{s-1} F_m(s) \rightarrow \Gamma(s) \zeta(s) \]
	locally uniformly for $\Re s>0$ (not just $\Re s>1$) as $m \rightarrow \infty$.
	This looks like a tall order. Crossing the pole at $s=1$ is not the problem. The difficulty
	is that unless $x$ is very close to 0, the expressions
	\[ \Delta_{m,k}(x) = \sum_{r=0}^k {k \choose r} (-1)^r p_m((r+1) x) \]
	involve values of $p_m$ at points well outside the interval $[0,1]$, where $p_m$ looks nothing
	like a negative exponential. Indeed $\Delta_{m,k}(x)$ is a divided difference of $p_m$ and consequently 
	equal to $(-x)^k p_m^{(k)}(u)$ for some $u$ between $x$ and $(k+1)x$ that is not easily specified.
	Since $p_m$ oscillates repeatedly on the interval $[0,m+1]$ it would seem that $\Delta_{m,k}(x)$ could
	be very large in size and of more or less random sign. So the following lemma comes as 
	something of a shock.
	
	\begin{lem}
		\label{magic}
		If $m$ is a non-negative integer and $p_m(x)=(1-x)(1-x/2)\ldots(1-x/m)$ then for each integer 
		$k$ and each $x \in [0,1]$
		\[ \Delta_{m,k}(x)=\sum_{r=0}^k {k \choose r} (-1)^r p_m((r+1) x) \geq 0.\]
	\end{lem}
	It is trivial to check that
	\[ \sum_{k=0}^m \Delta_{m,k}(x) = 1 \]
	for all $x$, so the lemma shows that for each $m$ the $\Delta_{m,k}$ form a partition of unity on
	$[0,1]$ and thus automatically controls the sizes of the $\Delta_{m,k}$ as well as their signs. 
	Once the lemma is established the convergence proof is fairly straightforward: this will be the content
	of Section \ref{conv} below.

	The obvious way to prove the Kronecker formula (\ref{Bern}) mentioned above is to use the expansion
	\[ \frac{y}{1-e^{-y}} = \sum_{k=0}^\infty \frac{(1-e^{-y})^{k}}{k+1} \]
	that already appeared in the integral formula for $\Gamma(s) \zeta(s)$. So it might be logically more reasonable to
	{\em define} the $F_m$ by using the formula
	\[ F_m(s)= \int_0^1 \left(\sum_{k=0}^m \frac{1}{k+1} \Delta_{m,k}(x) \right) x^{s-2} \, dx \] 
	and simply avoid mention of the Bernoulli numbers.
	However it seemed a little odd to define a rational function with known poles and residues as the
	analytic continuation of an integral.
	
	The convergence proof just alluded to relies on the fact that the $F_m$ are defined as sums
	which we can pass through integral signs. The point of the second section of the article will be to
	provide a bridge between the definition of the $F_m$ and their representation as characteristic polynomials: in other words to represent the $F_m$ as something more like a product than a sum.
	The main formula in Section \ref{bridge} is
	the following recurrence for the $F_m$: 
	
	\begin{lem}
		\label{recrel2}
		For each non-negative integer $m$
		\[ (s+m-1)F_m(s)=\frac{1}{m+1}+(m+1)\sum_{j=1}^m \frac{F_{m-j}(s)}{j(j+1)}. \]
	\end{lem}
	Thus
	\begin{eqnarray*}
		(s-1)F_0(s) & = & 1\\
		s F_1(s) & = & \frac{1}{2}+F_0(s) \\
		(s+1)F_2(s) & = & \frac{1}{3}+\frac{3}{2} F_1(s)+ \frac{3}{6} F_0(s)
	\end{eqnarray*}
	and so on.
	If we treat the first $m+1$ of these relations as a linear system for the
	values $F_0(s),F_1(s),\ldots,F_m(s)$ we can express the fact that $F_m(s)=0$ by the vanishing of a
	certain determinant. In Section \ref{spec} shall show that this determinant can be written as
	\[ \det(L_m+(1-s)U_m) \]
	where as mentioned earlier $L_m$ is the $(m+1)\times (m+1)$ Toeplitz matrix
	\[ \left( \begin{array}{cccccc} 1 & 0 & 0 & 0 & \ldots & 0 \\
	\frac{1}{2} & 1 & 0 & 0 & \ldots & 0 \\
	\frac{1}{3} & \frac{1}{2} & 1 & 0 & \ldots & 0 \\
	\frac{1}{4} & \frac{1}{3} & \frac{1}{2} & 1 &  & 0 \\
	\vdots & & & \ddots & \ddots & \vdots \\
	\frac{1}{m+1} & \frac{1}{m} & \frac{1}{m-1} & \ldots & \frac{1}{2} & 1 \end{array} \right) \]
	and $U_m$ is the matrix
	\[ \left( \begin{array}{cccccc} 0 & 1 & \frac{1}{2} & \frac{1}{3} & \ldots & \frac{1}{m} \\
	0 & 0 & \frac{1}{2} & \frac{1}{3} & \ldots & \frac{1}{m}\\
	0 & 0 &  0 & \frac{1}{3} & \ldots & \frac{1}{m} \\
	0 & 0 &  0 & 0 & & \vdots \\
	\vdots & \vdots & \vdots &  & \ddots & \frac{1}{m} \\
	0 & 0 & 0 & 0 & \ldots & 0 \end{array} \right). \]
	
	If we set $s=z/(z-1)$ the determinant becomes $\det(z L_m-(L_m+U_m))$ (apart from a factor of 
	$(1-z)^{m+1}$). The Riemann hypothesis would follow if this determinant vanishes only at points with $|z| \leq 1$ or equivalently that the matrix $I_{m+1}+L_m^{-1}.U_m$ has spectral radius at most 1, where $I_{m+1}$ is
	the $(m+1)\times (m+1)$ identity matrix. For small values of $m$ this is true. In the first version of this article I stated that there are good reasons to believe that the zeroes of the $F_m$ do leak across the critical line (and then come back again). Pace Nielsen quickly confirmed that when $m=643$ there is a zero to the right of the critical line. He also informed me that when $m=1127$ there is a zero with real part larger than 1 (which I had not expected). The Riemann hypothesis is equivalent to the statement that the spectral radius of $I_{m+1}+L_m^{-1}.U_m$ is at most $1+o(1)$ as $m \rightarrow \infty$.
	
	In Section \ref{est} of the article I shall explain why I think that the approximations $F_m(s)$ might be useful to estimate the size of $\zeta$. The main point is that whereas approximations to zeta that are sums of powers oscillate wildly all the way up the critical line, a polynomial of degree $m$ cannot oscillate too often.

	Whenever one has a new sequence of approximations to $\zeta$ it is natural to ask whether they can help to prove Diophantine properties (irrationality or transcendence) of values of the zeta function and most especially Euler's constant
	\[ \gamma = \lim_{s \rightarrow 1} \left( \zeta(s)-\frac{1}{s-1} \right). \]
	Since the approximations described here are rational functions (with integer coefficients) they do provide rational approximations to Euler's constant but for this particular ``value'' of zeta the approximations are not new.

	
	\section{The key lemma and convergence}
	\label{conv}
	
	In the introduction we defined, for each $m$, $p_m(x)=(1-x)(1-x/2)\ldots(1-x/m)$, and for each $k$
	\[ \Delta_{m,k}(x)=\sum_{r=0}^k {k \choose r} (-1)^r p_m((r+1) x). \]
	Note that the sum makes sense and is zero if $k<0$ or $k>m$. We introduced the function
	$F_m(s)$ as a rational function, 
	\[ F_m(s)=\sum_0^m \frac{a_{m,j} B_j}{s+j-1} \]
	and also mentioned the Kronecker formula 
	\[
	B_j = (-1)^j \sum_{k=0}^j \frac{1}{k+1} \sum_{r=0}^k {k \choose r} (-1)^r (r+1)^j. 
	\]
	It is a consequence of standard properties of the binomial coefficients that for all $j$ 
	between 0 and $m$, the sum on the right is unchanged if the upper limit is increased
	from $j$ to $m$. This implies that 
	\begin{eqnarray*} F_m(s)& = & 
		\sum_{j=0}^m (-1)^j \frac{a_{m,j}}{s+j-1} \sum_{k=0}^m \frac{1}{k+1} 
		\sum_{r=0}^k {k \choose r} (-1)^r (r+1)^j \\
		& = & \sum_{k=0}^m \frac{1}{k+1} \sum_{r=0}^k {k \choose r} (-1)^r 
		\sum_{j=0}^m (-1)^j \frac{a_{mj} (r+1)^j}{s+j-1} . \end{eqnarray*}
	For $\Re s>1$ the sum over $j$ can be written as
	\[ \sum_{j=0}^m (-1)^j a_{mj} (r+1)^j \int_0^1 x^{s+j-2} \, dx = \int_0^1 p_m((r+1) x) x^{s-2} \, dx \]  
	and so  for $\Re s>1$ we have
	\[ F_m(s)=\int_0^1 \sum_{k=0}^m \frac{\Delta_{m,k}(x)}{k+1} x^{s-2} \, dx. \]
	From now on we use $f_m(x)$ to denote $\sum_{k=0}^m \frac{\Delta_{m,k}(x)}{k+1}$.
	
	The aim of this section is to prove the following theorem 
	\begin{thm}
		\label{convthm}
		\[ h_m^{s-1} F_m(s) \rightarrow \Gamma(s) \zeta(s) \]
		\[ h_m^{s-1} (s-1) G_m(s) \rightarrow  \Gamma(s)  \]
		locally uniformly for $\Re s>0$ (with the obvious convention at $s=1$).
	\end{thm}
In view of the fact that $\Gamma$ has no zeros this theorem clearly implies the main convergence theorem of the article Theorem \ref{main}.

	It is easy to check that for $\Re s>1$
	\[ \Gamma(s) \zeta(s)= \int_0^\infty \frac{y}{1-e^{-y}} e^{-y} y^{s-2} \, dy. \]
	If we set $u=1-e^{-y}$ for $y>0$, then since $0<u<1$,
	\[
	\frac{y}{1-e^{-y}} =  \frac{-\log(1-u)}{u} 
	=  \sum_{k=0}^\infty \frac{u^{k}}{k+1} 
	= \sum_{k=0}^\infty \frac{(1-e^{-y})^{k}}{k+1}.
	\]
	Therefore 
	\[ \Gamma(s) \zeta(s)= \int_0^\infty \sum_{k=0}^\infty 
	\frac{1}{k+1} e^{-y} (1-e^{-y})^k y^{s-2} \, dy. \]
	
	Most of the effort in proving Theorem \ref{convthm}
	will go into showing that the truncated functions $f_m(x/h_m) {\bf 1}_{[0,h_m]}$ converge to 
	\[ \frac{x}{e^x-1}=\sum_{k=0}^\infty \frac{1}{k+1} e^{-x}(1-e^{-x})^k \]
	on $[0,\infty)$ with the convergence dominated by a negative exponential function.
	It is clear that for each fixed $x \geq 0$, $p_m(x/h_m) \rightarrow e^{-x}$ and hence that
	for each fixed $k$ and $x$
	\begin{equation}
	\Delta_{m,k}(x/h_m) \rightarrow e^{-x} (1-e^{-x})^k 
	\end{equation}
	as $m \rightarrow \infty$. We need to establish two types of dominance: one to confirm that the sum
	\[ f_m(x)=\sum_{k=0}^m \frac{\Delta_{m,k}(x)}{k+1} \]
	converges pointwise in $x$ to $\sum_{k=0}^\infty \frac{1}{k+1} e^{-x}(1-e^{-x})^k$ and one to check that this convergence is dominated on $[0,\infty)$.

	\setcounter{thm}{1}
	
	For almost every estimate we make it is essential to have the key lemma stated in the introduction:
	\begin{lem}
		If $m$ is a non-negative integer and $p_m(x)=(1-x)(1-x/2)\ldots(1-x/m)$ then for each integer 
		$k$ and each $x \in [0,1]$
		\[ \Delta_{m,k}(x)=\sum_{r=0}^k {k \choose r} (-1)^r p_m((r+1) x) \geq 0.\]
	\end{lem}
	We also need a simple
	property of the divided differences that depends only upon the fact that $p_m$ is a polynomial
	of degree at most $m$.
	
		\setcounter{thm}{4}
	
	\begin{lem}
		\label{upper}
		If $j$ is a non-negative integer then for every $x$,
		\[ \sum_{k=0}^m (k+1)(k+2)\ldots(k+j) \Delta_{m,k}(x) = j! p_m(-j x). \]
		In particular for $j=0$
		\[ \sum_{k=0}^m \Delta_{m,k}(x) =p_m(0)=1. \]
	\end{lem}
	
	\startproof
	We shall confirm that for any polynomial $q$ of degree at most $m$
	\[ \sum_{k=0}^m (k+1)(k+2)\ldots(k+j) \sum_{r=0}^k (-1)^r {k \choose r} q((r+1) x) = j! q(-j x) \]
	and in checking this we may assume that $x=1$.
	So our aim is to verify that for each such $q$
	\[ \sum_{k=0}^m (k+1)(k+2)\ldots(k+j) \sum_{r=0}^k (-1)^r {k \choose r} q(r+1) = j! q(-j). \]
	It suffices to check this for each polynomial of the form
	\[ q_n:t \mapsto (t-1)(t-2)\ldots(t-n) \]
	with $0 \leq n \leq m$. The internal sum vanishes if $q$ has degree less than $k$ and hence it vanishes for $q_n$ if $k>n$. It also vanishes if $k<n$ because of the form of $q_n$. The only remaining case is $k=n$ and in that case the internal sum has 
	value $(-1)^n n!$. So the double sum is
	\[ (n+1)(n+2)\ldots(n+j)(-1)^n n!=(-1)^n (n+j)!=j! q_n(-j). \]
	
	\endproof
	
	The proof of Lemma \ref{magic} involves the introduction of an additional
	parameter as follows. For each $v$ define
	\[ P_m(v,x)= m! p_m(x-v)=(v+1-x)(v+2-x)\ldots(v+m-x) \]
	and 
	\[ {\tilde \Delta}_{m,k}(v,x)=\sum_{r=0}^k  {k \choose r} (-1)^r P_m(v,(r+1) x). \]
	Observe that $P_m(0,x)=m! p_m(x)$ and ${\tilde \Delta}_{m,k}(0,x)=m! \Delta_{m,k}(x)$. 
	So Lemma \ref{magic} follows from:
	
	\begin{lem}
		If $m$ is a non-negative integer, $k$ is an integer, $v \geq 0$ and $0 \leq x \leq 1$
		\[ {\tilde \Delta}_{m,k}(v,x) \geq 0. \]
	\end{lem}
	
	\startproof
	We use induction on $m$. When $m=0$, ${\tilde \Delta}_{m,k}(v,x)$ is zero unless $k=0$ in which case
	it is 1. We claim that for $m>0$
	\[ {\tilde \Delta}_{m,k}(v,x) =(v+1-x){\tilde \Delta}_{m-1,k}(v+1,x)+
	k x {\tilde \Delta}_{m-1,k-1}(v+1-x,x) .\]
	Once this is established the inductive step is clear because we can assume that $k \geq 0$ and
	for the given range of $v$ and $x$, the number $v+1-x$ is also at least 0.
	
	Now for any $v$ and $x$
	\[ P_m(v,x)= (v+1-x) P_{m-1}(v+1,x) \]
	and so
	\[
	{\tilde \Delta}_{m,k}(v,x)  =  \sum_{r=0}^k {k \choose r} (-1)^r(v+1-(r+1)x) P_{m-1}(v+1,(r+1) x)\]
	\begin{eqnarray*}
		& = &  \sum_{r=0}^k  {k \choose r} (-1)^r (v+1-x) P_{m-1}(v+1,(r+1) x)  \\
		& & \hspace{0.5in}-\sum_{r=0}^k {k \choose r} (-1)^r r x P_{m-1}(v+1,(r+1) x)\\
		& = & (v+1-x)){\tilde \Delta}_{m-1,k}(v+1,x)  -  k \sum_{r=1}^k {k-1 \choose r-1} (-1)^r x P_{m-1}(v+1,(r+1) x)\\
		& = & (v+1-x)){\tilde \Delta}_{m-1,k}(v+1,x)  +k x \sum_{r=0}^{k-1} {k-1 \choose r} (-1)^r P_{m-1}(v+1,(r+2) x)\\
		& = &(v+1-x){\tilde \Delta}_{m-1,k}(v+1,x)+
		k x {\tilde \Delta}_{m-1,k-1}(v+1-x,x)
	\end{eqnarray*}
	where the last step follows from the fact that for all $m$, $v$ and $x$,
	\[ P_{m-1}(v+1,(r+2)x)=P_{m-1}(v+1-x,(r+1)x). \]
	
	\endproof

	By combining Lemmas \ref{magic} and \ref{upper} we can immediately make some estimates for the $  \Delta_{m,k}(x)$
	that will give us part of the dominance we need to get convergence.
	
	\begin{lem}
		\label{sumdom}
		For each $m \geq 0$, each $k \geq 0$ and each $x \in [0,h_m]$. 
		\[ \Delta_{m,k}(x/h_m) \leq \frac{1}{(k+1)} e^{x}\]
		and for each $m \geq 1$
		\[ f_m(x/h_m)-p_m(x/h_m)= \sum_{k=1}^m \frac{1}{k+1} \Delta_{m,k}(x/h_m) \leq e^{x}-1. \]
	\end{lem}
	
	\startproof For the first one we apply Lemma \ref{upper} with $j=1$ and use the positivity of the  
	$\Delta_{m,k}$ to deduce that for each $k$
	\[ \Delta_{m,k}(x) \leq \frac{1}{k+1} p_m(-x) = \frac{1}{k+1} (1+x)(1+x/2)\ldots(1+x/m) \leq 
	\frac{1}{(k+1)} e^{h_m x}. \]
	For the second one we observe that 
	\[ \sum_{k=1}^m \frac{1}{k+1} \Delta_{m,k}(x)\leq  \sum_{k=0}^m k \Delta_{m,k}(x) = p_m(-x)-p_m(0)
	\leq e^{h_m x}-1. \]
	
	\endproof
	
	As already remarked it is clear that for each fixed $x \geq 0$, $p_m(x/h_m) \rightarrow e^{-x}$ and hence that
	for each fixed $k$ and $x$
	\begin{equation}
	\label{xkconv} \Delta_{m,k}(x/h_m) \rightarrow e^{-x} (1-e^{-x})^k 
	\end{equation}
	as $m \rightarrow \infty$. From Lemma \ref{sumdom} we have
	\[ \Delta_{m,k}(x/h_m) \leq \frac{1}{(k+1)} e^{x}\]
	so the convergence in (\ref{xkconv}) is dominated (on the space of non-negative integers 
	with counting measure)
	by a sequence summable against $(1/(k+1))$. Hence for each $x>0$
	\[ f_m(x/h_m)=\sum_{k=0}^m \frac{1}{k+1} \Delta_{m,k}(x/h_m) \rightarrow \frac{x}{e^x-1}. \]
	
	We have that 
	\[ h_m^{s-1} F_m(s) = \int_0^{h_m} f_m(x/h_m) x^{s-2} \, dx. \]
	In order to use dominated convergence on $[0,\infty)$ we need an estimate for $f_m$ which we prove
	by introducing another extra parameter.
	For each $m$ and $p \geq 0$ define
	\[ {\tilde f}_m(p,x)=\sum_{k=0}^m \frac{1}{k+1+p} \Delta_{m,k}(x) \] 
	and observe that $f_m(x)={\tilde f}_m(0,x)$.

	\begin{lem}
		\label{pest}
		For each $m$ and $p$
		\[ {\tilde f}_m(p,x)= {\tilde f}_{m-1}(p,x)+\frac{p x}{m} {\tilde f}_{m-1}(p,x)-
		\frac{(p+1)x}{m} {\tilde f}_{m-1}(p+1,x). \]
	\end{lem}

	\startproof
	As long as $p>-1$ we have
	\begin{eqnarray*}
		{\tilde f}_m(p,x) & = & \sum_{k=0}^m \Delta_{m,k}(x) \int_0^{\infty} \frac{1}{(1+u)^{k+2+p}} \, du \\
		& = & \int_0^{\infty} K_{m,x}(u) \frac{1}{(1+u)^{p+1}} \, du 
	\end{eqnarray*}
	where
	\[ K_{m,x}(u) = \sum_{k=0}^m \frac{1}{(1+u)^{k+1}} \Delta_{m,k}(x) . \]
	This function is holomorphic on the plane apart from its pole at $-1$. Its derivatives at 0
	are successively
	\[ \sum_{k=0}^m \Delta_{m,k}(x) \]
	\[ -\sum_{k=0}^m (k+1) \Delta_{m,k}(x) \]
	\[ \sum_{k=0}^m (k+1)(k+2) \Delta_{m,k}(x) \]
	and so on and therefore by Lemma \ref{upper}  its power series expansion at 0 is
	\[ \sum_{j=0}^\infty (-1)^j p_m(-j x) u^j. \]
	
	Therefore, for $|u|<1$
	\begin{eqnarray*} K_{m,x}(u) & = &\sum_{j=0}^\infty (-1)^j p_m(-j x) u^j \\
		& = & \sum_{j=0}^\infty (-1)^j p_{m-1}(-j x)\left( 1+\frac{j x}{m} \right) u^j \\
		& = & K_{m-1,x}(u)+\frac{x}{m} \sum_{j=0}^\infty (-1)^j p_{m-1}(-j x) \, j u^j \\
		& = & K_{m-1,x}(u)+ \frac{x}{m}  u K_{m-1,x}'(u). 
	\end{eqnarray*}
	The outermost identity continues analytically so it holds for all $u>-1$. Therefore
	\begin{eqnarray*}
		{\tilde f}_m(p,x) & = & \int_0^{\infty} K_{m,x}(u) \frac{1}{(1+u)^{p+1}} \, du \\
		& = & {\tilde f}_{m-1}(p,x)+ \frac{x}{m}
		\int_0^{\infty} u K_{m-1,x}'(u) \frac{1}{(1+u)^{p+1}} \, du \\
		& = & {\tilde f}_{m-1}(p,x)- \frac{x}{m}
		\int_0^{\infty} K_{m-1,x}(u) \frac{d}{du} \frac{u}{(1+u)^{p+1}} \, du \\
		& = & {\tilde f}_{m-1}(p,x)+ \frac{x}{m} \int_0^{\infty} K_{m-1,x}(u) 
		\left( \frac{p}{(1+u)^{p+1}} -\frac{p+1}{(1+u)^{p+2}} \right)\, du \\
		& = & {\tilde f}_{m-1}(p,x)+ \frac{p}{m} x {\tilde f}_{m-1}(p,x)-
		\frac{p+1}{m} x {\tilde f}_{m-1}(p+1,x). 
	\end{eqnarray*}
	
	\endproof

	Now we can estimate $f_m(x)$ as follows. Using the key lemma we have an inequality
	\[ {\tilde f}_m(1,x)=\sum_{k=0}^m \frac{1}{k+2} \Delta_{m,k}(x) \geq 
	\frac{1}{2} \sum_{k=0}^m \frac{1}{k+1} \Delta_{m,k}(x) = \frac{1}{2} f_m(x) \]
	provided $0 \leq x \leq 1$.
	Then by Lemma (\ref{pest})
	\[ f_m(x)={\tilde f}_m(0,x)= {\tilde f}_{m-1}(0,x)-\frac{x}{m}  {\tilde f}_{m-1}(1,x) \leq
	f_{m-1}(x) \left(1-\frac{x}{2 m} \right) \leq e^{-x/(2 m)} f_{m-1}(x). \]
	So by induction
	\[ f_m(x) \leq e^{-h_m x/2} \]
	and we get the negative exponential dominance
	\[ f_m(x/h_m) \leq e^{-x/2} \]
	on the range of integration $[0,h_m]$.
	This suffices to guarantee that for $\Re s>1$ 
	\[ h_m^{s-1} F_m(s) = \int_0^{h_m} f_m(x/h_m) x^{s-2} \, dx \rightarrow \Gamma(s) \zeta(s). \]
	
	We wish to cross the pole and so we need to modify the integrand.
	For $\Re s>1$ we have
	\begin{eqnarray*}
		 \Gamma(s) \left( \zeta(s)-\frac{1}{s-1} \right) & = & \int_0^{\infty} \frac{x}{e^x-1} x^{s-2} \, dx -
	\int_0^{\infty} e^{-x} x^{s-2} \, dx \\
	& = &
	\int_0^{\infty} \left( \frac{x}{e^x-1}-\frac{1}{e^x} \right) x^{s-2} \, dx. \end{eqnarray*}
	The last integrand behaves like $x$ near 0 so the integral converges locally uniformly for
	$\Re s>0$ and represents the holomorphic function $\Gamma(s) \left( \zeta(s)-\frac{1}{s-1} \right) $
	on this larger region.
	
	As in the introduction we set 
	\[ G_m(s) = \sum_{j=0}^m (-1)^j \frac{a_{m,j}}{s+j-1} \]
	and observe that for $\Re s>1$
	\[ G_m(s)= \int_0^1 p_m(x) x^{s-2} \, dx. \]
	So on this half-plane
	\[ h_m^{s-1} G_m(s) = \int_0^{h_m} p_m(x/h_m) x^{s-2} \, dx. \]
	Now $p_m(x/h_m) \rightarrow e^{-x}$ for each fixed $x$ and $0 \leq p_m(x/h_m) \leq e^{-x}$
	as long as $0 \leq x \leq h_m$ so for $\Re s>1$
	\[ \int_0^{h_m} p_m(x/h_m) x^{s-2} \, dx \rightarrow \Gamma(s) \frac{1}{s-1}. \]
	Therefore, still only for $\Re s>1$, we have
	\[  h_m^{s-1} (F_m(s) - G_m(s)) = 
	\int_0^{h_m} \left( f_m(x/h_m)-p_m(x/h_m) \right) x^{s-2} \, dx \rightarrow 
	\Gamma(s) \left( \zeta(s)-\frac{1}{s-1} \right) . \]
	The integrand is dominated as $x \rightarrow \infty$ by $e^{-x/2}$ but also, as $x \rightarrow 0$
	by $ e^{x}-1$
	owing to Lemma \ref{sumdom}. Moreover, $F_m$ and $G_m$ both have residue 1 at $s=1$
	so the difference is holomorphic for $\Re s>0$. So the integral represents $h_m^{s-1} (F_m(s) - G_m(s))$
	on the larger region and also converges to  $\Gamma(s) \left (\zeta(s)-\frac{1}{s-1} \right) $ on this region.
	
	To complete the proof of Theorem \ref{convthm} it suffices to prove the (easy) second assertion of the theorem:
	\[ (s-1) h_m^{s-1} G_m(s) \rightarrow \Gamma(s) \]
	locally uniformly for $\Re s>0$. We have that for $\Re s>1$,
	\[(s-1) h_m^{s-1} G_m(s) = (s-1) \int_0^{h_m} p_m(x/h_m) x^{s-2} \, dx =
	-\int_0^{h_m} 1/h_m p_m'(x/h_m) x^{s-1} \, dx  \]
	because $p_m(1)=0$. The latter integral converges as long as $\Re s>0$ so it represents 
	$(s-1) h_m^{s-1} G_m(s)$ on the larger region. It suffices to show that
	\[ -1/h_m p_m'(x/h_m) {\bf 1}_{[0,h_m]} \rightarrow e^{-x} \]
	for $x>0$ and that the convergence is dominated by a negative exponential.
	
	Observe that $p_m$ is decreasing on $[0,1]$ so $-p_m'(x)$ is positive for $0 \leq x < 1$.
	On the other hand
	\[ -p_m'(x) = p_m(x) \sum_{j=1}^m \frac{1}{j-x} \leq \prod_{j=2}^m (1-x/j)+
	p_m(x) \sum_{j=2}^m \frac{1}{j-1} \leq e^{-h_m x}(1+h_m).  \]
	Thus for $0 \leq x \leq h_m$
	\[ -1/h_m p_m'(x/h_m)  \leq e^{-x} (1/h_m+1) \leq 2 e^{-x} \]
	giving the required dominance.
	Also
	\[ -1/h_m p_m'(x/h_m)  = \frac{1}{h_m} \prod_{j=2}^m (1-x/(j h_m))+
	p_m(x/h_m) \frac{1}{h_m} \sum_{j=2}^m \frac{1}{j-x/h_m}. \]
	The first term is at most $\frac{1}{h_m} e^{-(h_m-1)x/h_m}$ and tends to 0 while the
	second term behaves like $e^{-x} (h_m-1)/h_m$ and tends to $e^{-x}$.
	This establishes Theorem \ref{convthm}.
	
	Finally we have that the ratios
	\[ \frac{F_m(s)}{(s-1)G_m(s)} \]
	converge locally uniformly to $\zeta(s)$ for $\Re s>0$.
	My guess is that they do so on the entire complex plane.
	
	\section{The bridge from sum to determinant}
	\label{bridge}
	\setcounter{thm}{2}
	The purpose of this section is to establish the recurrence
	\begin{lem}
	For each non-negative integer $m$
	\[ (s+m-1)F_m(s)=\frac{1}{m+1}+(m+1)\sum_{j=1}^m \frac{F_{m-j}(s)}{j(j+1)} \]
\end{lem}
\setcounter{thm}{8}
	 which will enable us
	to express the numerator of $F_m$ as a determinant. A similar recurrence holds for the functions $G_m$: if $m \geq 1$
	\[ (s+m-1)G_m(s)=(m+1)\sum_{j=1}^m \frac{G_{m-j}(s)}{j(j+1)}. \]
	A small modification of the proof below actually yields this as well.
	
	We begin with a simple remark.
	\begin{lem}
		For each non-negative integer $m$
		we have 
		\[ \lim_{s \rightarrow \infty} s F_m(s) = \frac{1}{m+1}. \]
	\end{lem}
	
	\startproof 
	\[ \lim_{s \rightarrow \infty} s F_m(s) = f_m(1). \]
	All the $\Delta_{m,k}$ vanish at $x=1$ apart from $\Delta_{m,m}$ since they involve only values of
	$p_m$ at the integers $1,2,\ldots,m$. By Lemma \ref{upper} the $\Delta_{m,k}$ add up to 1 so
	$\Delta_{m,m}(1)=1$. So 
	\[ f_m(1)=\sum_{k=0}^m \frac{\Delta_{m,k}(1)}{k+1}=\frac{1}{m+1}. \] 
	\endproof
	
	Now for the proof of Lemma \ref{recrel2}.

	\startproof
	The two sides have the same limits at infinity so it suffices to check that they have the 
	same residues at each of the points $1,0,-1,\ldots,2-m$. At $s=1-r$ the residue on the
	left is $(m-r) a_{m,r} B_r$ while the residue on the right is
	\[ (m+1) \sum_{j=1}^m \frac{a_{m-j,r} B_r}{j(j+1)}. \]
	It thus suffices to check that for each $0 \leq r \leq m-1$,
	\[ (m-r)a_{m,r} = (m+1) \sum_{j=1}^m \frac{a_{m-j,r}}{j(j+1)}. \]
	(The case $r=m$ is also obvious.)
	Multiplying by $x^r$ and summing over $0 \leq r \leq m$ it is enough to check that
	\[ m \, p_m(x)- x p_m'(x)=(m+1) \sum_{j=1}^m \frac{p_{m-j}(x)}{j(j+1)}. \]
	Both sides are polynomials of degree $m-1$ so we need only check the values at
	$x=1,2,\ldots,m$.
	Let $k$ be one of these integers. $p_m(k)=0$ and it is easy to check that
	\[ -k \, p_m'(k)= (-1)^{k-1} {m \choose k}^{-1}. \]
	On the other hand
	\[ (m+1) \sum_{j=1}^m \frac{p_{m-j}(k)}{j(j+1)}= 
	(m+1) \sum_{j=m-k+1}^m \frac{1}{j(j+1)} \prod_{i=1}^{m-j} \left(1 - \frac{k}{i} \right) \]
	because $p_{m-j}$ vanishes at $k$ if $m-j \geq k$.
	The latter expression is
	\begin{eqnarray*}
		(m+1) \sum_{j=m-k+1}^m \frac{(-1)^{m-j}}{j(j+1)} {k-1 \choose m-j} & = &
		(m+1) \sum_{r=0}^{k-1} \frac{(-1)^{r}}{(m-r)(m-r+1)} {k-1 \choose r} \end{eqnarray*}
	\begin{eqnarray*}
		& = & (m+1) \sum_{r=0}^{k-1} {k-1 \choose r} (-1)^{r} \int_0^1 (x^{m-r-1}-x^{m-r}) \, dx \\
		& = & -(m+1) \int_0^1 x^m \left(1-\frac{1}{x} \right)^k \, dx \\
		& = & (-1)^{k-1} (m+1) \int_0^1 x^{m-k} \left(1-x \right)^k \, dx 
	\end{eqnarray*}
	and the beta integral gives the appropriate reciprocal of the binomial coefficient. 
	\endproof 
	
	The recurrence relation given by Lemma \ref{recrel2} describes a dynamical system for the sequence $F_m(s)$. Numerically this system appears to evolve very slowly and indeed the convergence of the sequence is very slow. This makes the approximations useless for effective calculation of the value $\zeta(s)$ but suggests that it might be possible to track the dynamical system: it is almost a continuous-time system.

	\section{The spectrum}
	\label{spec}
	
	From the previous section we have that for each $m$
	\begin{equation} \label{recrel} (s+m-1)F_m(s)=\frac{1}{m+1}+(m+1)\sum_{j=1}^m \frac{F_{m-j}(s)}{j(j+1)}. \end{equation}
	The first $m+1$ of these relations give us the linear system
	\[ \left( \begin{array}{cccccc}
	s-1 & 0 & 0 & 0 & \ldots & 0 \\
	-\frac{2}{2} & s & 0 & 0 & \ldots & 0 \\
	-\frac{3}{6} & -\frac{3}{2} & s+1 & 0 & \ldots & 0 \\
	\vdots & & & \ddots & & \vdots \\
	\vdots & & & & \ddots & \vdots \\
	-\frac{m+1}{m(m+1)} & \ldots & \ldots & -\frac{m+1}{6} & -\frac{m+1}{2} & s+m-1 
	\end{array} \right) 
	\left( \begin{array}{c} F_0(s) \\ F_1(s) \\ \vdots \\  \\ F_m(s) \end{array} \right) =
	\left( \begin{array}{c} 1 \\ \frac{1}{2} \\ \vdots \\ \\ \frac{1}{m+1} \end{array} \right). \] 
	So $F_m(s)$ can be written as the ratio of two determinants. The denominator is the determinant
	of the $(m+1) \times (m+1)$ matrix on the left, namely $(s-1)s(s+1) \ldots (s+m-1)$. The numerator 
	is the determinant of the matrix obtained by replacing the last column of the original 
	matrix with the vector 
	$(1,1/2,\ldots, 1/(m+1))$. It will be more convenient to move this vector to the first column of the
	matrix thus introducing a factor of $(-1)^m$ into the determinant, and to change the sign of all the 
	other columns, thus removing the factor again. Then the numerator can be written

	\[ \det \left( \begin{array}{ccccccc}
	1 & 1-s & 0 & 0 & 0 & \ldots & 0 \\
	\frac{1}{2} & \frac{2}{2} & -s & 0 & 0 & \ldots & 0 \\
	\frac{1}{3} & \frac{3}{6} & \frac{3}{2} & -s-1 & 0 & \ldots & 0 \\
	\vdots & & & & \ddots & & \vdots \\
	\vdots & & & & & \ddots & \vdots \\
	\frac{1}{m} & \frac{m}{(m-1) m} & \ldots & & \frac{m}{6} & \frac{m}{2} & -s-m+2 \\ 
	\frac{1}{m+1} & \frac{m+1}{m(m+1)} & \ldots & & & \frac{m+1}{6} & \frac{m+1}{2} 
	\end{array} \right). \]
	Regard the columns as labelled $0,1,\ldots,m$. We leave the zero and 1 columns unchanged as
	$(1,1/2,1/3,\ldots,1/(m+1))$ and $(1-s,1,1/2,\ldots,1/m)$ respectively. We add the 1 column to the
	2 column to get 
	\[ \left( (1-s), (1-s), 2, 1,2/3, \ldots, 2/(m-1) \right). \] 
	We add the (new) 2 column to the 3 column and we get
	\[ \left( (1-s), (1-s),(1-s), 3,3/2,3/3, \ldots, 3/(m-2) \right).\] 
	
	Continue in this way
	and after all the additions divide the 2 column by 2, the 3 column by 3 and so on to get
	\begin{equation} 
	\label{Fmatrixform} F_m(s)=\frac{m! \det(L_m+(1-s)U_m)}{(s-1)s \ldots (s+m-1)} \end{equation}
	where
	\begin{equation} \label{log} L_m = \left( \begin{array}{cccccc} 1 & 0 & 0 & 0 & \ldots & 0 \\
	\frac{1}{2} & 1 & 0 & 0 & \ldots & 0 \\
	\frac{1}{3} & \frac{1}{2} & 1 & 0 & \ldots & 0 \\
	\frac{1}{4} & \frac{1}{3} & \frac{1}{2} & 1 &  & 0 \\
	\vdots & & & \ddots & \ddots & \vdots \\
	\frac{1}{m+1} & \frac{1}{m} & \frac{1}{m-1} & \ldots & \frac{1}{2} & 1 \end{array} \right)  
	\end{equation}
	and 
	\begin{equation} \label{int}
	 U_m = \left( \begin{array}{cccccc} 0 & 1 & \frac{1}{2} & \frac{1}{3} & \ldots & \frac{1}{m} \\
	0 & 0 & \frac{1}{2} & \frac{1}{3} & \ldots & \frac{1}{m}\\
	0 & 0 &  0 & \frac{1}{3} & \ldots & \frac{1}{m} \\
	0 & 0 &  0 & 0 & & \vdots \\
	\vdots & \vdots & \vdots &  & \ddots & \frac{1}{m} \\
	0 & 0 & 0 & 0 & \ldots & 0 \end{array} \right). \end{equation}
	Notice that $F_m$ does not have poles at $-2,-4,\ldots$ because the corresponding
	Bernoulli numbers vanish. So the determinant in the numerator picks up the trivial zeros of
	zeta at these numbers: indeed the factor $s+2$ appears as soon as $m \geq 3$, the factor $s+4$ as soon as
	$m \geq 5$ and so on.
	
	The zeroes of $F_m$ are thus related in a simple way to the spectrum of $L_m^{-1} U_m$. 
	This formulation for $F_m$ is perhaps the most elegant one in terms of
	a determinant but it is interesting to express $F_m$
	in a form in which the problem looks more like a conventional spectral problem.
	
	To begin with we multiply column $j$ by $j$ for each $j \geq 1$ but leave the zero column 
	unchanged. This includes the $m!$ factor appearing in (\ref{Fmatrixform}) into the determinant.
	We now subtract the $m-1$ row from the last, the $m-2$ row from the $m-1$ and so on
	to produce the matrix
	\[  \left( \begin{array}{cccccc} 1 & 0 & 0 & 0 &\ldots & 0 \\
	-\frac{1}{2} & 1 & 0 & 0 & \ldots & 0 \\
	-\frac{1}{6} & -\frac{1}{2} & 2 & 0 & \ldots & 0 \\
	-\frac{1}{12} & -\frac{1}{6} & -\frac{1}{2} & 3  & & 0 \\
	\vdots & & & \ddots & \ddots & \vdots \\
	\frac{-1}{m(m+1)} & \frac{-1}{(m-1)m} & \frac{-2}{(m-2)(m-1)} & \ldots &-\frac{1}{2} & m 
	\end{array} \right) +
	(1-s)\left( \begin{array}{cccccc} 0 & 1 & 1 & 1 & \ldots & 1 \\
	0 & -1 & 0 & 0 & \ldots & 0 \\
	0 & 0 & -1 & 0 & & \vdots \\
	\vdots & \vdots &  & \ddots & \ddots & 0 \\
	0 & 0  & \ldots & 0 & -1 & 0 \\
	0 & 0 & 0  & \ldots & 0 & -1 \end{array} \right). \]
	We now add all the rows below the top one, to the top one, so that the second matrix now has a
	zero top row. The first matrix now has top row
	\[ \left( \frac{1}{m+1}, \frac{1}{m}, \frac{2}{m-1}, \ldots, \frac{m}{1} \right). \]
	Since the variable $s$ now appears on the diagonal in all places except the first, our aim is to
	reduce the dimension by one so as to create a characteristic polynomial proper. We add multiples
	$(m+1)/2$, $(m+1)/6$ and so on of the top row to the successive rows below. Since this eliminates
	the first column below the first row, the determinant is now the top left entry $1/(m+1)$
	multiplied by the determinant of the remaining square. So we get
	\[ \frac{\det(T_m+R_m+(s-1) I_m)}{m+1} \]
	where $T_m$ and $R_m$ are as follows.
	\[ T_m= \left( \begin{array}{ccccc} 
	1 & 0 & 0 & \ldots & 0 \\
	-\frac{1}{2} & 2 & 0 & \ldots & 0 \\
	-\frac{1}{6} & -\frac{2}{2} & 3 & & \vdots \\
	\vdots & \vdots & & \ddots & 0 \\
	\frac{-1}{(m-1)m} & \frac{-2}{(m-2)(m-1)} & \ldots & -\frac{m-1}{2} & m 
	\end{array} \right) \]
	which is lower triangular with entry $j$ in the $jj$ diagonal place and entry 
	\[ -j/((i-j)(i-j+1)) \]
	in the $ij$ place, if $i>j$. $R_m$ is the rank one matrix given by
	\[ R_m = \left( \frac{(m+1)j}{i(i+1)(m-j+1)} \right)_{ij}. \]
	We have
	\[ F_m(s)=\frac{1}{m+1} \frac{\det(T_m+R_m+(s-1) I_m)}{\prod_{j=0}^m (s+j-1)}. \]
	We are thus interested in the spectrum of the matrix $T_m+R_m$ where $T_m$ is a certain lower triangular matrix and $R_m$ has rank 1. It will be seen that this formulation is actually somewhat closer to the recurrence relation (\ref{recrel}). 
	
	If we set
	\[ A_m = T_m+R_m \; \; \; \mbox{ and } \; \; \; B_m = T_m+R_m-I_m \]
	then the determinant in question is 
	\[ \det (s A_m-(s-1)B_m). \]
	The complex numbers to the right of the critical line are those for which $|s|>|s-1|$ so a natural way to tackle the spectral problem would be to try to find a norm on $\C^m$ with the property that for every $z \in C^m$
	\[ \| A_m z \| \geq \|B_m z \|. \]
	The obvious choice would be an Hilbertian norm. So we look for a positive definite matrix $H$ for which
	\[ A_m^* H A_m -B_m^* H B_m \]
	is also positive definite or alternatively one for which
	\[ A_m H A_m^* -B_m H B_m^* \]
	is positive definite.
	
	The form of the matrix $T_m$ makes this very tempting. If we ignore the rank one matrix $R_m$ then we can certainly compute the spectrum of $T_m$ since it is lower triangular. However there is a natural choice of norm which shows that the spectrum is to the left of the critical line and therefore provides a more robust argument that one could try to perturb. If $H$ is the diagonal matrix with entries $1,1/2,1/3,\ldots,1/m$ on the diagonal then 
	\[  T_m H T_m^*- (T_m-I_m) H (T_m-I_m)^* = T_m H+H T_m^*-H\]
	is the matrix 
	\[ \left( \begin{array}{ccccc} 
	1 & -\frac{1}{2} & -\frac{1}{6} & \ldots & \frac{-1}{(m-1)m} \\
	-\frac{1}{2} & \frac{3}{2} & -\frac{1}{2} & \ldots &  \frac{-1}{(m-2)(m-1)}\\
	-\frac{1}{6} & -\frac{1}{2} & \frac{5}{3} & & \vdots \\
	\vdots & \vdots & \ddots & \ddots & -\frac{1}{2} \\
	\frac{-1}{(m-1)m} & \frac{-1}{(m-2)(m-1)} & \ldots & -\frac{1}{2} & 2-\frac{1}{m} 
	\end{array} \right). \]
	This is obviously positive definite because it has negative off diagonal entries and row sums that are positive
	because of the familiar telescoping sum
	\[ \sum_{n=1}^m \frac{1}{n(n+1)} = 1-\frac{1}{m+1}. \]
	
	As remarked earlier there are good reasons to think that the zeroes of the $F_m$ do leak through the critical line so that the best one can hope for is to find matrices $H_m$ with
	\[ (1+\epsilon_m) A_m^* H_m A_m -B_m^* H_m B_m \]
	positive definite and $ \epsilon_m \rightarrow 0$.
	
	Once one is in possession of the matrices $T_m$ and $R_m$ one could confirm that they yield the $F_m$ in a ``direct'' way. Diagonalise $T_m$ as $PDP^{-1}$ and check that when you apply $P^{-1}$ and $P$ to the components of $R_m$ you recover the Bernoulli and Stirling numbers. Such an argument would be a bit harsh on the reader since there would be little motivation for introducing these particular matrices. More importantly the dynamical system described by the recurrence relation (\ref{recrel}) is of interest in itself.
	
	\section{Estimating the size of $\zeta$}
	\label{est}
	Numerical evidence indicates that the function $f_m(x/h_m)$ differs from $x/(e^x-1)$ by only about $1/(h_m)^2$ at any point of
	$[0,h_m]$ and so we expect the ratio 
	\[ \frac{F_m(s)}{(s-1)G_m(s)} \]
	to provide a good approximation to $\zeta$ at $s=1/2+it$ as long as $\Gamma(s)$ is as large as $1/(h_m)^2$. This happens if $|t|$ is at most a bit less than
	$ \frac{2}{\pi} \log \log m$.
	In fact, numerical evidence and rough calculations indicate that the ratio is not too far from $\zeta$ for $t$ all the way up to a multiple of $\log m$. At the same time there are good reasons to think that $F_m(s)$ does not oscillate significantly for $t$ larger than $\log m$. So we have the tantalising possibility that the two regions overlap: the $t < \log m$ region where $F_m$ tells us about $\zeta$ and the $t > \log m$ region where $F_m$ is smooth enough to be estimated.
	
	This discussion suggests that one should look at the asymptotic expansion for $F_m(s)$ which starts off
	\[ F_m(s)= \frac{1}{(m+1)(s-1)}+ \frac{c_m}{(m+1)(s-1)s}+ \cdots \]
	where the coefficient $c_m$ grows logarithmically with $m$, the next coefficient like $(\log m)^2$ and so on. However my feeling is that the more promising approach is the ``usual'' one: to look at an integral (say)
	\[ F_m(s)-G_m(s) = \int_0^1 (f_m(x)-p_m(x)) x^{s-2} \, dx \]
	and move the contour into the region where $x^{s-2}$ is very small if $s=1/2+i t$ with $t$ large.
	
	For the genuine $\zeta$ integral
	\[ \int_0^\infty \frac{x}{1-e^{-x}} e^{-x} x^{s-2} \, dx \]
	this approach is hopeless because the contour is forced up against the imaginary axis and hence picks up the poles of
	\[ z \mapsto \frac{z}{1-e^{-z}}. \] 
	Being a polynomial, $z \mapsto f_m(z)$ has no poles so the issue does not arise. The problem is to estimate $f_m$ off the real line.
	
	\section{The connection with the Connes, Berry-Keating operator}
	\label{Connes}
	
	The articles by Connes  \cite{con} and Berry and Keating \cite{BK} each describe an operator on an infinite-dimensional space whose spectrum shares properties of the zeros of the zeta function. The operators are formally the same but the spaces on which they are considered are different. Connes showed that all Riemann zeros that lie on the critical line correspond to eigenvalues of his operator but was not able to check this for zeros off the line (if they exist). Berry and Keating showed that the mean density of their eigenvalues matches that of the Riemann zeros.
	
	In Connes' incarnation, the operator can be built from a multiplication operator and an integral operator acting on an infinite-dimensional function space as explained in the article \cite{lach}. From this one can see that the finite-dimensional operators considered in this article are sections of the Connes operator as follows.
	
	The Toeplitz matrix $L_m$ given in equation (\ref{log}) can be thought of as acting on 
	polynomials $a_0+a_1 x+a_2 x^2+\cdots+a_m x^m$ rather than sequences $(a_0,\ldots,a_m)$. It does so by multiplication by the partial sum
	\[ \sum_0^{m} \frac{x^j}{j+1} \]
	of the series for $\frac{- \log(1-x)}{x} $ (followed by truncation back to a polynomial of degree $m$).
	In this context the upper triangular matrix $U_m$ in (\ref{int}) maps the constant function to 0 and for each
	$k \geq 1$ the monomial $x^k$ to the sum
	\[ \frac{1+x+\cdots+x^{k-1}}{k} = \frac{1-x^k}{k(1-x)}. \]
	Thus for any polynomial $q$ of degree $m$ the image is
	\[ \frac{1}{1-x} \int_x^1 \frac{q(t)-q(0)}{t} \, dt. \]
	
	The finite-dimensional operators described here have several advantages. 
	\begin{itemize}
	\item The determinants yield approximations to the zeta function itself (not just eigenvalues that correspond to zeros).
	
	\item Having finite-dimensional operators means that one need not worry about the space on which the matrices act: although of course we would like to find the right norm in order to prove things about the eigenvalues.
	
	\item These approximations provably do converge to zeta and so in the limit, pick up all the Riemann zeros.
	
	\end{itemize}
	
	To create finite-dimensional sections of an operator by restricting and projecting onto polynomials is normally an extremely unstable thing to do unless one is working with very special normed spaces (such as $L_2$ of the disc). The fact that it works here is rather remarkable and may perhaps indicate that when trying to find a norm in order to check the spectrum one should start with something like $\ell_2^m$.  
	
	\section*{Acknowledgements}
	I am extremely grateful to David Preiss for his advice during this work. My thanks also to Terry Tao and Peter Sarnak who made very helpful suggestions concerning the presentation.

	\noindent Keith Ball \\
	Mathematics Institute \\
	University of Warwick \\
	Coventry CV4 7AL
	
	\noindent k.m.ball@warwick.ac.uk

\end{document}